\documentclass[11pt]{amsart}
\usepackage[mathscr]{eucal}
\usepackage{amsfonts}
\usepackage{amsmath}
\usepackage{amsthm}
\usepackage{amssymb}
\usepackage{latexsym}

\theoremstyle{plain}
\newtheorem{definition}[equation]{Definition}
\newtheorem{conjecture}[equation]{Conjecture}
\newtheorem{corollary}[equation]{Corollary}
\newtheorem{lemma}[equation]{Lemma}
\newtheorem{prop}[equation]{Proposition}
\newtheorem{theorem}[equation]{Theorem}

\theoremstyle{definition}
\newtheorem{remark}[equation]{Remark}

\numberwithin{equation}{subsection}

\newcommand{\C}{\mathbb C}

\newcommand{\F}{\mathcal F}
\newcommand{\GG}{{\mathbf{\Gamma}_{n}}}

\newcommand{\hh}{\mathsf{H}}

\newcommand{\ot}{\otimes}
\newcommand{\PP}{\mathbb P}
\newcommand{\pa}{\partial}
\newcommand{\Rep}{{\mathrm{Rep}}}

\newcommand{\ZZ}{\mathbb Z}

\def\A{{\mathsf{A}}}
\def\R{{\mathsf{R}}}
\def\E{{\mathsf{E}}}
\def\k{{\Bbbk}}
\def\QQ{{\overline{Q}}}
\newcommand{\otot}{{\otimes\cdots\otimes}}
\newcommand{\ve}{\varepsilon}
\newcommand{\lf}{\lfloor}
\newcommand{\rf}{\rfloor}
\newcommand{\into}{\,\,\hookrightarrow\,\,}
\newcommand{\vi}{${\sf {(i)}}\;$}
\newcommand{\vii}{${\sf {(ii)}}\;$}

\begin{document}
\title{Generalized double affine Hecke algebras of higher rank}

\author{Pavel Etingof}
\address{Department of Mathematics, Massachusetts Institute of
Technology, Cambridge, MA 02139, U.S.A.}
\email{etingof@math.mit.edu}

\author{Wee Liang Gan}
\address{Department of Mathematics, Massachusetts Institute of
Technology, Cambridge, MA 02139, U.S.A.}
\email{wlgan@math.mit.edu}

\author{Alexei Oblomkov}
\address{Department of Mathematics, Massachusetts Institute of
Technology, Cambridge, MA 02139, U.S.A.}
\email{oblomkov@math.mit.edu}

\maketitle{}

\section{Introduction}

Double affine Hecke algebras (DAHA) appeared in the work of Cherednik
\cite{Ch}, as a tool to prove Macdonald's conjectures; since that
time they have been in the center of attention of many
representation theorists. In particular, in \cite{Sa}, Sahi
extended them to root systems of type $C^\vee C_n$, and used
this extension to establish Macdonald's conjectures for
Koornwinder polynomials.

DAHA have a rich algebraic structure, which relates them to
algebraic geometry and the theory of integrable systems.
For instance, it is shown in \cite{Ob1}
that Cherednik's DAHA of type $A_{n-1}$ is a quantization of
the relativistic Calogero-Moser space (the space of states for the
Ruijsenaars-Schneider integrable system). Also, it is shown
in \cite{Ob2} that Sahi's DAHA of rank 1 is a quantization
of a generic cubic surface with three lines forming a triangle
removed.

Motivated by this, P.E., A.O., and Eric Rains
introduced generalized DAHA (GDAHA) of rank 1, attached to
any star-shaped affine Dynkin diagram, i.e.
$\widetilde D_4,\widetilde E_6,\widetilde E_7,\widetilde E_8$,
\cite{EOR}. It was shown in \cite{EOR} that the PBW theorem
holds for these algebras, and that they provide quantizations
of del Pezzo surfaces (with a singular genus one curve removed).
In the case of $\widetilde D_4$, the GDAHA is the same
as the Sahi algebra (of rank 1), so one recovers the results of
\cite{Ob2}.

Later, it was pointed out in \cite{ER} that the definition
of GDAHA of rank 1 makes sense and the PBW property remains true
for any star-shaped graph $D$ which is not a Dynkin diagram of finite
type. Such GDAHA are flat deformations of group algebras of
polygonal Fuchsian groups acting on the Euclidean (in the affine
case) or the Lobachevsky plane.

The main goal of this paper is to introduce and begin to study
GDAHA of higher rank $n>1$, attached to any star-shaped
graph that is not a finite Dynkin diagram.
These algebras are deformations of the semidirect products
of the symmetric group $S_n$ with the $n$-th tensor power
of the rank 1 GDAHA. Like the original DAHA, GDAHA are quotients
 of group algebras of appropriate braid groups,
and reduce to the Sahi algebras (of rank $n$) in the case of
the graph $\widetilde D_4$.

This paper is organized as follows.
In Section 2, we study degenerate, or rational GDAHA of higher rank.
These algebras are not really new, as they are ``spherical''
subalgebras of the algebras introduced by W.L.G. and V. Ginzburg
in \cite{GG} (Def. 1.2.3), associated to the idempotent of the branching
vertex; in the affine case they are also ``spherical''
subalgebras of the wreath product symplectic reflection algebras
introduced in \cite{EG}. However, we give a new presentation of
rational GDAHA by generators and relations, which is a higher rank
generalization of Theorem 1 in \cite{Me}  (see also \cite{MOV}).
Using this presentation, we give a parametrization
of irreducible representations of rational GDAHA
for affine $D$ when the quantum parameter $\hbar$ vanishes.
This parametrization is by the space of solutions
of a certain additive Deligne-Simpson problem, which turns out to
be a smooth algebraic variety of dimension $2n$.
Since for affine $D$ and generic parameters
GDAHA are Morita equivalent to symplectic reflection
algebras, this parametrization is not really new,
and essentially coincides with the parametrization of
representations by generalized
Calogero-Moser spaces, i.e. by quiver-theoretical data given in
\cite{EG}, Theorem 11.16; however, our new presentation is
somewhat simpler.

In Section 3, we define GDAHA and prove a formal PBW theorem for
them. The proof is based on the fact that formal GDAHA are
a special case of Hecke algebras of orbifolds introduced in
\cite{E}, for which the formal PBW theorem holds in a very
general situation. On the other hand, the algebraic PBW theorem
(i.e. the freeness of the algebra as a module over the ring
of coefficients) unfortunately remains a conjecture.

In Section 4, we introduce the Knizhnik-Zamolodchikov connection
with coefficients in the degenerate GDAHA, and use it to
construct the monodromy functor from the category of finite
dimensional representations for the degenerate GDAHA to
that for the nondegenerate one; this gives a large supply
of finite dimensional representations of GDAHA in the affine
case, by applying the monodromy functor to the representations
from \cite{EM,M,Ga}. This connection
also allows us to construct a Riemann-Hilbert homomorphism
between completions of the nondegenerate and the degenerate
GDAHA, which gives another (more elementary) proof of the formal
PBW theorem. We note however that although the formal PBW theorem
for GDAHA is a purely algebraic statement, both proofs we
give are based on the Riemann-Hilbert correspondence and
therefore use complex analysis; we don't know a purely
algebraic proof.

In Section 5, we study GDAHA in the case of affine $D$.
Unfortunately, outside of type $\widetilde D_4$,
we are unable to establish any of the important
properties of GDAHA (proved for usual DAHA of type A in
\cite{Ob1}), and they are stated as conjectures.
Basically, we expect that GDAHA are quantizations
of spaces of Calogero-Moser type, which are
(topologically trivial) deformations of Hilbert schemes
of affine del Pezzo surfaces described above.

The main result of Section 5 is
the construction of the parametrization
of irreducible representations of GDAHA for $q=1$
by points of generalized relativistic Calogero-Moser spaces,
which are defined as spaces of solutions of a
certain multiplicative Deligne-Simpson problem.
This parametrization is the non-degenerate analog of
the parametrization in the degenerate case
constructed in Section 2, and we conjecture it to be a
bijection. With respect to these two parametrizations,
the monodromy functor induces the usual Riemann-Hilbert map
between solutions of the additive and multiplicative
Deligne-Simpson problems (see e.g. \cite{CB2}).

We note that if $m$ is the number of ``legs'' of $D$, then
the monodromy maps discussed above depend
on a choice of $m$ points $\alpha_1,...,\alpha_m$ on $\Bbb C\Bbb P^1$,
modulo fractional linear transformations.
Thus, in the case $\widetilde {D}_4$, there is an essential
parameter -- the cross-ratio $\kappa$ (while in the other affine
cases there is no such parameter). If one keeps the value of the monodromy
fixed and varies $\kappa$, one gets a flow on the space
of solutions of the additive Deligne-Simpson problem (which is
$2n$-dimensional). In the case $n=1$ this is the Painlev\'e VI
flow (see e.g. \cite{EOR}, Section 7), so for $n>1$ this flow
should be regarded as a higher rank analogue of Painlev\'e VI.
Note that this flow for $n>1$ has an additional parameter (so it
has 5 rather than 4 parameters); it would be interesting to
study this flow in more detail, for instance write it down
explicitly.

We expect that for any non-Dynkin graph there exist multiplicative analogs of
the Gan-Ginzburg algebras (\cite{GG}, Definition 1.2.3),
which are higher rank generalizations
of multiplicative preprojective algebras introduced in
\cite{CBS}. Moreover, we expect that GDAHA are
``spherical'' subalgebras of the multiplicative Gan-Ginzburg
algebras corresponding to the idempotent of the branching
vertex. In the rank 1 case, this is shown by Crawley-Boevey and
Shaw in the appendix to \cite{EOR}.

{\bf Acknowledgments.} The work of P.E. and A.O. was partially
supported by the NSF grant DMS-9988796 the CRDF grant
RM1-2545-MO-03. The work of W.L.G. was partially supported by the
NSF grant  DMS-0401509. The authors want to thank Jasper Stokman
for the interest to the work and suggestions which improved the
exposition.

\section{Rational GDAHA}

\subsection{Conventions}

Throughout the paper, given an affine algebraic
variety $X$ over $\Bbb C$, a point $x\in X$, and a
module $M$ over $\Bbb C[X]$, we will write $M(x)$ for the
specialization of $M$ at $x$.

Let $D$ be a star-like graph, i.e. a tree with one $m$-valent
vertex, called the node or the branching vertex, and the rest of the vertices
2- and 1-valent (which form $m$ ``legs'' growing from the node).
We will always assume that $D$ is {\bf not} a finite Dynkin diagram.
We label the legs of $D$ by numbers $1,...,m$.
Let $d_k$ be the number of vertices in the $k$-th leg
of $D$, including the node. Let $I$ be the set of vertices of
$D$, and let $i_0\in I$ be the node.
Let $i_1(k), \ldots, i_{d_k-1}(k)$ be the vertices of the
$k$-th leg of $D$ enumerated from the node.

Let $\gamma=(\gamma_{kj})$, $k=1,...,m$, $j=1,...,d_k$,
be a collection of variables, and $\nu$ an additional variable.
It is easy to show that there exist unique $\mu_i=\mu_i(\gamma)$,
$i\in I$, and $\xi_k=\xi_k(\gamma)$, $k=1,...,m$ such that
$\xi_1+...+\xi_m=0$ and
\begin{equation}\label{mugamma}
\gamma_{kj}=\sum_{p=1}^{j-1}\mu_{i_p}
+\frac{\mu_{i_0}}{m}+\xi_k
\end{equation}
for all $j,k$.

\subsection{The definition of the rational GDAHA}

\begin{definition}
The rational (or degenerate) generalized DAHA of rank $n$ attached
to $D$ is the algebra $B_n$ generated over $\C[\gamma,\nu]$ by
elements $Y_{i,k}$ (where $i=1,\ldots,n$; $k=1,\ldots,m$) and the symmetric
group $S_n$, with the following defining relations:
for any $i,j,h \in [1,n]$ with $i\neq j$, and $k, l\in [1,m]$,

\begin{gather*}
s_{ij}Y_{i,k} = Y_{j,k}s_{ij}, \\
s_{ij}Y_{h,k} = Y_{h,k}s_{ij} \quad\textrm{if $h\neq i,j$},\\
\prod_{j=1}^{d_k}(Y_{i,k}-\gamma_{kj})=0,\\
Y_{i,1}+Y_{i,2}+\cdots+Y_{i,m}=\nu\sum_{j\neq i}s_{ij},\\
[Y_{i,k},Y_{j,k}] = \nu(Y_{i,k}-Y_{j,k})s_{ij},\\
[Y_{i,k},Y_{j,l}] = 0,\ k\ne l,
\end{gather*}
where $s_{ij}\in S_n$ denotes the transposition
$i\leftrightarrow j$.
\end{definition}

\begin{remark} It is obvious that
the algebra $B_n(\gamma,\nu)$ does not change
(up to an isomorphism) under the transformations
of parameters $\gamma_{kj}\to \gamma_{kj}+\sigma_k$, where
$\sigma_k\in \Bbb C$, $\sigma_1+...+\sigma_m=0$
(the required isomorphism is given by $Y_{i,k}\to Y_{i,k}+\sigma_k$).
So the essential parameters of $B_n(\gamma,\nu)$ are $\mu_i$
and $\nu$, and there are $m$ ``redundant''
parameters $\xi_k$. However, it is convenient to keep the redundant
parameters to simplify the presentation. A similar remark
applies to GDAHA defined in Section 3 below.
\end{remark}

\subsection{GDAHA and the Gan-Ginzburg algebras}

We now recall some definitions from \cite{GG}.
Let $\k$ be a commutative ring.
Let $Q$ be a quiver, and denote by
$I$ the set of vertices of $Q$.
The double $\QQ$ of $Q$ is the quiver obtained from $Q$ by
adding a reverse edge $\stackrel{a^*}{j\to i}$
for each edge $\stackrel{a}{i\to j}$ in $Q$.
If $\stackrel{a}{i\to j}$ is an edge
in $\QQ$, we call $t(a):=i$ its tail, and
$h(a):=j$ its head.

Let $R:=\bigoplus_{i\in I} \k$,
and $E$ be the free $\k$-module with
basis formed by the set of edges $\{a\in \QQ\}$.
Thus, $E$ is naturally a $R$-bimodule and $E=\bigoplus_{i,j\in I}
E_{i,j}$, where $E_{i,j}$ is spanned by the edges $a\in\QQ$ with
$h(a)=i$ and $t(a)=j$. The path algebra of $\QQ$
is $\k\QQ := T_R E = \bigoplus_{n\geq 0} T^n_R E$, where
$T^n_R E = E\ot_R \cdots \ot_R E$ is the $n$-fold
tensor product.
The trivial path for the vertex $i$ is denoted by
$e_i$, an idempotent in $R$.

Let $n$ be a positive integer. Let $\R := R^{\ot n}$.
For any $\ell\in [1,n]$, define the $\R$-bimodules
$$ \E_{\ell}:= R^{\ot (\ell-1)}\ot E\ot R^{\ot (n-\ell)}
\qquad \mathrm{and}\qquad \E :=
\bigoplus_{1\leq\ell\leq n} \E_\ell\,.$$
The natural inclusion $\E_\ell \into
R^{\ot (\ell-1)}\ot T_R E\ot R^{\ot (n-\ell)} $
induces a canonical identification
$T_\R\E_\ell = R^{\ot (\ell-1)}\ot T_R E\ot R^{\ot (n-\ell)}$.
Given two elements $\ve\in\E_\ell$ and $\ve'\in\E_m$ of the form
\begin{equation} \label{eqnve}
\ve= e_{i_1}\ot e_{i_2}\otot a\otot h(b) \otot e_{i_n}\,,
\end{equation}
\begin{equation} \label{eqnve'}
\ve'= e_{i_1}\ot e_{i_2}\otot t(a)\otot b\otot e_{i_n}\,,
\end{equation}
where $\ell\neq m$, $a,b\in \QQ$ and $i_1, \ldots, i_n \in I$,
we define
\begin{eqnarray}
\lf\ve,\ve'\rf &:=&
(e_{i_1}\otot a\otot h(b) \otot e_{i_n})
(e_{i_1}\otot t(a)\otot b \otot e_{i_n}) \nonumber\\ & &
-(e_{i_1}\otot h(a)\otot b \otot e_{i_n})
(e_{i_1}\otot a\otot t(b) \otot e_{i_n}) . \nonumber
\end{eqnarray}
Note that $\lf\ve,\ve'\rf$ is an element in $T^2_\R \E$.

\begin{definition}[\cite{GG}, Defn. 1.2.3] \label{dea}
For any $\mu=(\mu_i)_{i\in I}$, where $\mu_i\in \k$,
and $\nu\in \k$, define the algebra
$\A_n=\A_{n,\mu,\nu}$ to be the quotient of $T_{\R}\E\rtimes\k[S_n]$
by the following relations.
\begin{itemize}
\item[\vi]
For any $i_1,\ldots, i_n\in I$ and $\ell\in [1,n]$:
\begin{gather*}
e_{i_1}\otot \left(
\sum_{\{a\in Q\,|\, h(a)=i_\ell\}}
a\cdot a^*- \sum_{\{a\in Q\,|\, t(a)=i_\ell\}} a^*\cdot a
-\mu_{i_\ell}e_{i_\ell} \right)
\otot e_{i_n}  \\
= \nu \sum_{\{ j\neq\ell \,|\, i_j=i_\ell\}}
(e_{i_1}\otot e_{i_\ell} \otot e_{i_n})s_{j\ell}.
\end{gather*}
\item[\vii]
For any $\ve, \ve'$ of the form (\ref{eqnve})--(\ref{eqnve'}):
\[  \lf\ve,\ve'\rf = \left\{ \begin{array}{ll}
\nu(e_{i_1}\otot h(a)\otot t(a)\otot e_{i_n})s_{\ell m}
& \textrm{if $b\in Q$,\ $a = b^*$} ,\\
- \nu(e_{i_1}\otot h(a)\otot t(a)\otot e_{i_n})s_{\ell m}
& \textrm{if $a\in Q$,\ $b = a^*$} ,\\
0 & \textrm{else}\,.  \end{array} \right.  \]
\end{itemize}
\end{definition}

Now, given a star-like graph $D$,
we let $Q=Q(D)$ be the quiver obtained from $D$
by assigning an orientation to the edges of $D$
so that they look away from the node.
Let $e_i$ be the idempotent in $\C^I$
that corresponds to the vertex $i$.
Denote by $\mathsf{A}_n$ the algebra over $\k:=\Bbb C[\mu,\nu]$
associated to $Q$ defined Definition \ref{dea}.

Let $\xi=(\xi_1,...,\xi_m)$ be a set of variables such that
$\xi_1+...+\xi_m=0$.

\begin{prop} \label{prop20}
There is a natural isomorphism $\varphi: B_n \simeq
e_{i_0}^{\ot n}\mathsf{A}_ne_{i_0}^{\ot n}\otimes \Bbb C[\xi]$.
\end{prop}

\begin{proof} In the case when $n=1$ and $D$ is affine,
the proof is given in \cite{Me}, see also \cite{MOV} and
\cite[Prop. 7.2]{EOR}. In general, the proof is analogous.

Namely, note that the algebra $B_n$ has a natural filtration
defined by the condition $\deg(Y_{i,k})=2$; similarly
$\mathsf{A}_n$ has a filtration defined by giving the edges of
$\overline{Q}$ degree 1.
Let $h_k$ be the edge of $\overline{Q}$ that starts at
$i_0$ and goes along the $k$-th leg of $\overline{Q}$.
Let $h_k^*$ be the edge of $\overline{Q}$ opposite to $h_k$.
Then we can define a filtration preserving homomorphism
$\varphi: B_n \to
e_{i_0}^{\ot n}\mathsf{A}_ne_{i_0}^{\ot n}\otimes \Bbb C[\xi]$
by the formula
$$
\varphi(Y_{i,k})
=e_{i_0}^{\ot i-1}\ot (h_k^*h_k+(\xi_k+\frac{\mu_{i_0}}{m})e_{i_0})
\ot e_{i_0}^{\ot n-i}
$$
for all $i,k$.
It is clear that ${\rm gr}(B_n)$ is a quotient of
$B_n(0,0)[\gamma,\nu]$, and
it is shown in \cite{GG} that ${\rm gr}\mathsf{A}_n=
\mathsf{A}_n(0,0)[\mu,\nu]$. Also,
by Theorem 1 in \cite{Me},
the specialization of $\varphi$ at $(0,0)$ is an isomorphism.
Thus ${\rm gr}\varphi$ is an isomorphism and hence
$\varphi$ is an isomorphism, and the proposition is proved.
\end{proof}

Proposition \ref{prop20} and the results of \cite{GG} imply the
following corollary.

\begin{corollary}\label{fre}
The natural homomorphism
$B_n(0,0)[\gamma,\nu]\to {\rm gr}B_n$
is an isomorphism, and thus $B_n$ is a
free $\C[\gamma,\nu]$-module.
\end{corollary}

Assume now that $D$ is affine.
Let $\hh(\GG)$ be the symplectic reflection algebra
associated (as in \cite{EG}) to the wreath product group
 $\GG = S_n\ltimes \Gamma^n$, where
$\Gamma$ is the finite subgroup of $SL(2,\C)$ corresponding
to $D$. Let $V_{i}$ be the representation of $\Gamma$
corresponding to the vertex $i\in I$ under the McKay correspondence, and
let $e_i\in \C[\Gamma]$ be a primitive idempotent of this
representation. Let $p=\sum_i e_i$.
Then $p^{\otimes n}\in \C[\GG]$

It is shown in \cite{GG} that
$p^{\otimes n}\hh(\GG)p^{\otimes
n}=\mathsf{A}_n$. This fact and
Proposition \ref{prop20} imply the
following.

\begin{corollary} \label{cor20}
There is a natural isomorphism
$\psi: B_n \to e_{i_0}^{\otimes n}\hh(\GG)e_{i_0}^{\otimes
n}\otimes \Bbb C[\xi]$.
\end{corollary}

\subsection{The degenerate cyclotomic Hecke algebra}\label{dcha}

In this subsection we define the degenerate version of
the Ariki-Koike cyclotomic Hecke algebra.

Let $\ell$ be a positive integer, and
$\lambda=(\lambda_1,...,\lambda_\ell)$ be a collection of
variables.

\begin{definition} The degenerate cyclotomic Hecke algebra
${\bold B}_{n,\ell}$ is the algebra over
$\C[\lambda,\nu]$,
generated by $S_n$ and additional generators $Y_1,...,Y_n$,
with defining relations
\begin{gather*}
s_{ij}Y_{i} = Y_{j}s_{ij}, \\
s_{ij}Y_{h} = Y_{h}s_{ij} \quad\textrm{if $h\neq i,j$},\\
(Y_{i}-\lambda_1)\cdots
(Y_{i} -\lambda_\ell)=0,\\
[Y_{i},Y_{j}] =
\nu(Y_{i}-Y_{j})s_{ij}.
\end{gather*}
\end{definition}

Note that ${\bold B}_{n,\ell}(\lambda,0)=\C[S_n]\ltimes
R_\lambda^{\otimes n}$, where $R_\lambda=\Bbb C[Y]/(\prod
(Y-\lambda_j))$.

For each $k=1,...,m$, we have a homomorphism
$\eta_k: {\bold B}_{n,d_k}\to B_n$, such that
$\eta_k(\lambda)=\gamma_k$, where
$(\gamma_k)_j:=\gamma_{kj}$, and $\eta_k(\nu)=\nu$.
This homomorphism is given by
the formulas $\eta_k(Y_i)=Y_{i,k}$, $\eta_k(s_{ij})=s_{ij}$.
It is easy to check that the specialization of $\eta_k$
at $(0,0)$ is injective.
By Corollary \ref{fre}, this implies that
the natural map ${\bold B}_{n,\ell}(0,0)[\lambda,\nu]\to
{\rm gr}{\bold B}_{n,\ell}$ is an isomorphism, and hence
${\bold B}_{n,\ell}$ is a free module over
$\C[\lambda,\nu]$ of rank $n!\ell^n$. Thus
${\bold B}_{n,\ell}(\lambda,\nu)$ is an algebra of dimension $n!\ell^n$
for all $\lambda,\nu$, which is semisimple for generic values of
parameters.

The algebra ${\bold B}_{n,\ell}(\gamma,\nu)$ has a 1-dimensional
representation $\chi$ given by the formula
$\chi(s_{ij})=1$, $\chi(Y_i)=\gamma_\ell$.
We call this representation the
trivial representation.
For generic parameters $\lambda,\nu$,
the representation $\chi$ defines an
idempotent in ${\bold B}_{n,\ell}(\lambda,\nu)$.
We will denote this idempotent by ${\rm e}$.

For any representation
$V\in {\rm Rep} {\bold B}_{n,\ell}(\gamma,\nu)$, we
denote by $V^{{\bold B}_{n,\ell}(\gamma,\nu)}$ the space of
homomorphisms of representations $\chi\to V$.
Obviously, this space can be naturally regarded
as a subspace of $V$. In the generic case, this subspace is
equal to ${\rm e}V$.

The algebra ${\bold B}_{n,\ell}(\gamma,\nu)$ contains an obvious
subalgebra ${\bold B}_{n-1,\ell}(\gamma,\nu)$, generated by $Y_i$ and
$s_{ij}$ with $i,j<n$. For any representation $V$
of ${\bold B}_{n,\ell}(\gamma,\nu)$, denote by $V'$ the space
$V^{{\bold B}_{n-1,\ell}(\gamma,\nu)}$.

Consider the element $x:=Y_n-\nu\sum_{j=1}^{n-1}s_{nj}$.
It is easy to check that $x$ commutes with
 ${\bold B}_{n-1,\ell}(\gamma,\nu)$, hence $x$ preserves the space
$V'$.

Let $T$ be the $n$-by-$n$ matrix such that
$T_{ij}=1-\delta_{ij}$. Note that
$T+1$ has rank $1$.

\begin{lemma}\label{eigen} Let $V$ be the regular representation of
${\bold B}_{n,\ell}(\gamma,\nu)$. Then
for generic parameters, the space $V'$ has dimension
$n\ell$, and $x|_{V'}$ is conjugate to
$$
(\lambda_{\ell} {\rm Id}_n-\nu T)
\oplus\mathop{{\rm diag}}(\lambda_{1},\dots,\lambda_{\ell-1})\otimes
{\rm Id}_n.
$$
\end{lemma}

\begin{proof}
Let $\pi$ be the idempotent in ${\bold
B}_{n-1,\ell}(\gamma,\nu)\subset {\bold B}_{n,\ell}(\gamma,\nu)$
corresponding to the character $\chi$.
Let $P_j(x):=\prod_{p\ne
j}\frac{x-\lambda_p}{\lambda_j-\lambda_p}$. It is easy to check
(by considering the case $\nu=0$)
that the elements $v_{ij}:=\pi P_j(Y_n)s_{ni}$
(where we agree that $s_{nn}=1$) form a basis of
$V'$. Let us compute the action of $x$ in this basis.
If $i\ne n$, we have
\begin{gather*}
xv_{ij}=\pi xP_j(Y_n)s_{ni}=
\lambda_j v_{ij}-\nu\sum_{p\ne n}\pi s_{np}P_j(Y_n)s_{ni}=\\
\lambda_j v_{ij}-\nu\sum_{p\ne n}\pi P_j(Y_p)s_{np}s_{ni}=
\lambda_j v_{ij}-\nu\delta_{j\ell}\sum_{p\ne n}\pi
s_{np}s_{ni}=\\
\lambda_j v_{ij}-\nu\delta_{j\ell}\pi (\sum_{p\ne n,i}
s_{ip}s_{np}+1)=\\
\lambda_j v_{ij}-\nu\delta_{j\ell}\pi \sum_{p\ne
i}\sum_{q=1}^\ell P_q(Y_n)s_{np}=
\lambda_j v_{ij}-\nu\delta_{j\ell}\sum_{p\ne
i}\sum_{q=1}^\ell v_{pq}.
\end{gather*}
If $i=n$, we have the same result:
\begin{gather*}
xv_{nj}=\pi xP_j(Y_n)=
\lambda_j v_{nj}-\nu\sum_{p\ne n}\pi s_{np}P_j(Y_n)=\\
\lambda_j v_{nj}-\nu\sum_{p\ne n}\pi P_j(Y_p)s_{np}=
\lambda_j v_{nj}-\nu\delta_{j\ell}\sum_{p\ne n}\pi
s_{np}=\\
\lambda_j v_{nj}-\nu\delta_{j\ell}\pi \sum_{p\ne
n}\sum_{q=1}^\ell P_q(Y_n)s_{np}=
\lambda_j v_{nj}-\nu\delta_{j\ell}\sum_{p\ne
n}\sum_{q=1}^\ell v_{pq}.
\end{gather*}
This implies the required statement.
\end{proof}

\subsection{The affine case}\label{affi}

Consider now the affine case, i.e.
$D=\widetilde D_4,\widetilde E_6,\widetilde E_7,\widetilde E_8$.
Then $m=3,4$ and the numbers $d_k$ are the following
(up to ordering): $(2,2,2,2),(3,3,3),(2,4,4)$, and $(2,3,6)$,
respectively. We let $\delta_i$ be the coordinates of the basic
imaginary root $\delta$ of $D$ in the basis of simple roots.
Also, let $0\in I$ be the vertex corresponding to the trivial
representation of $\Gamma$ under the McKay correspondence.

Let us assume that $\ell:=d_m$ is the largest of the $d_k$.
In this case $\ell$ is divisible by $d_k$ for all $k$. Set
$$
\hbar=\hbar(\gamma):=\ell\sum_{k,j}\frac{\gamma_{kj}}{d_k}.
$$

The main properties of $B_n$ in the affine case
are summarized in the following theorem.

\begin{theorem}\label{mapro} (i) The Gelfand-Kirillov dimension
of $B_n(\gamma,\nu)$ is $2n$.

(ii) The algebra $B_n(\gamma,\nu)$ is PI if and only
if $\hbar=0$. If $\hbar=0$, this algebra is PI of degree
$n!\ell^n$.

(iii) If $\hbar=0$, then $B_n(\gamma,\nu)$ is finitely
generated over its center $Z(B_n(\gamma,\nu))$.
Moreover, for generic $\gamma,\nu$ 
(with $\hbar=0$) the map $Z(B_n(\gamma,\nu))\to {\rm e}B_n(\gamma,\nu){\rm e}$
given by $z\mapsto z{\rm e}$ is an isomorphism. In particular,
${\rm e}B_n(\gamma,\nu){\rm e}$ is a commutative algebra.

(iv) If $\hbar=0$ and otherwise $(\gamma,\nu)$ are generic
then $B_n(\gamma,\nu)$
is an Azumaya algebra, and ${\mathcal R}_{n,\gamma,\nu}:
={\rm Spec}(Z(B_n(\gamma,\nu)))$
is a smooth affine algebraic variety of dimension $2n$. In this
case, every irreducible representation of $B_n(\gamma,\mu)$
restricts (via the map $\eta_m$) to the regular representation of
$\bold B_{n,\ell}(\gamma_m,\nu)$.
\end{theorem}

\begin{proof}
(i) By Corollary \ref{cor20}, 
the associated graded algebra of $B_n$ under its natural
filtration is 
$$
e_{i_0}^{\otimes n}(\Bbb C[\bold \Gamma_n]\ltimes 
\Bbb C[x_1,...,x_n,y_1,...,y_n])e_{i_0}^{\otimes
n}.
$$
This implies that the Gelfand-Kirillov dimension of
$B_n$ is $2n$.

(ii) If $\hbar=0$ and otherwise $\gamma,\nu$ are generic, 
then by Theorem 16.1 of \cite{EG} and Corollary \ref{cor20}, 
$B_n(\gamma,\nu)$ is an Azumaya algebra,
whose fibers are matrix algebras of size $n!\ell^n$. 
Thus, $B_n(\gamma,\nu)$ is PI of this degree. Hence it is PI of
degree $\le n!\ell^n$ for any $\gamma,\nu$ with $\hbar=0$. 
But the associated graded of $B_n$ is clearly PI of degree
exactly $n!\ell^n$, so the statement follows. 

On the other hand, let $\hbar\ne 0$. The algebra $B_n(\gamma,\nu)$ induces 
a Poisson bracket on the center $Z_0=\Bbb
C[x_1,...,x_n,y_1,...,y_n]^\GG$ 
of the algebra $B_n(0,0)$. 
It follows from \cite{EG}, Section 2, that this Poisson bracket 
is the one induced by a symplectic form on $\Bbb C^{2n}$. 
The Poisson center of $Z_0$ under this bracket consists 
only of scalars. This implies that the center
$Z(B_n(\gamma,\nu))$ is trivial, and hence $B_n(\gamma,\nu)$ 
is not PI. 

(iii) By Corollary \ref{cor20}, the center $Z(B_n)$ of $B_n$ coincides
with the center of the symplectic reflection algebra 
$\hh(\GG)$. It is proved in \cite{EG} that 
$\hh(\GG)$ is finite over its center, so the first statement of
(iii) follows. The rest follows from the proof of (ii) (the
Azumaya property of $B_n$). 

(iv) The first two statements follow from the Azumaya property
of $B_n$ and Section 11 of \cite{EG}. To prove the last
statement, note that since $\bold B_{n,\ell}(\gamma_m,\nu)$
is a semisimple algebra for generic parameters, 
it is sufficient to prove the statement for $\nu=0$ 
and generic representations. In this case, the result
follows easily from the rank 1 case, see \cite{CBH}. 
\end{proof}

\subsection{Representations of $B_n(\gamma,\nu)$ for $\hbar=0$.}\label{rep0}

Assume that $\hbar=0$ and otherwise $(\gamma,\nu)$ are generic.
Theorem 11.16 of \cite{EG} furnishes an isomorphism
of algebraic varieties $\Phi_{EG}: {\mathcal R}_{n,\gamma,\nu}\to {\bold
M}_{n,\gamma,\nu}$ of ${\mathcal R}_{n,\gamma,\nu}$
onto a certain explicitly described variety ${\bold
M}_{n,\gamma,\nu}$
(the Calogero-Moser space attached to $D$), which is a deformation of the
Hilbert scheme of the desingularization of the Kleinian
singularity $\Bbb C^2/\Gamma$. By the definition,
${\bold M}_{n,\gamma,\nu}$ is the variety
of isomorphism classes of
representations of the doubled quiver $\overline{Q}$
with dimension vector $n\delta$
such that $\sum_{a\in Q} [a,a^*]=\sum \mu_ie_i-\nu Te_0$
(here $\mu_i$ are related to $\gamma_{kj}$
by formula (\ref{mugamma})).

It is convenient for us to give a slightly different description
of the variety ${\bold M}_{n,\gamma,\nu}$.

\begin{definition} Define ${\mathcal  M}_{n,\gamma,\nu}$
to be the variety of conjugacy classes of $m$-tuples $(x_1,\dots,x_m)\in
{\frak {gl}}_{n\ell}(\C)^{m}$ satisfying the following equations:
\begin{gather}\label{producl}
x_1+x_2+...+x_m=0,\\
\label{gam_12..}
x_k\sim \mathop{{\rm diag}}(\gamma_{k1},\dots,\gamma_{kd_k})\otimes
Id_{n\ell/d_k},\quad k=1,\dots,m-1,\\
\label{x_m}
x_m\sim (\gamma_{m\ell} {\rm Id}_n-\nu T)
\oplus\mathop{{\rm diag}}(\gamma_{m1},\dots,\gamma_{m,\ell-1})\otimes
{\rm Id}_n.
\end{gather}
That is, ${\mathcal  M}_{n,\gamma,\nu}$ is the categorical quotient
of the variety $\widetilde {\mathcal M}_{n,\gamma,\nu}$
of $m$-tuples as above by the action of the group
$PGL_{n\ell}(\Bbb C)$.

Here, if $x,y\in {\frak gl}_N(\C)$,
we use the notation $x\sim y$ to say that $x$ and
$y$ are in the same conjugacy class.
\end{definition}

\begin{remark}
Thus, ${\mathcal  M}_{n,\gamma,\nu}$ is defined as the variety
of solutions of an appropriate additive Deligne-Simpson problem.
\end{remark}

\begin{prop}\label{smocl}
For generic parameters ${\mathcal  M}_{n,\gamma,\nu}$ is a smooth
variety, of dimension $2n$.
\end{prop}

\begin{proof}
The proof is standard and analogous to the proof of
Proposition \ref{smo} below.
\end{proof}

\begin{prop}\label{beta} There exists a regular map
$\beta: {\bold M}_{n,\gamma,\nu}\to {\mathcal M}_{n,\gamma,\nu}$
which sends a representation of $\overline{Q}$
from ${\bold M}_{n,\gamma,\nu}$ to the collection of operators \linebreak
$h_k^*h_k+\xi_k+\frac{\mu_{i_0}}{m}$, where $h_k,h_k^*$ are
defined in the proof of Proposition \ref{prop20}.
This map is an isomorphism.
\end{prop}

\begin{proof} The proof is based on the following lemma from
linear algebra, due to Crawley-Boevey.

Let $\Lambda_i\in \Bbb C$
for $i=1,...,N$,
such that $\Lambda_i+...+\Lambda_j\ne 0$ for any $1\le i\le j\le N$.
Let $V_i$, $i=0,...,N+1$, $N\ge 0$,
be finite dimensional complex vector spaces
of dimensions $D_i$, $D_{i-1}<D_i$. Let ${\mathcal O}$ be a
conjugacy class in ${\frak {gl}}(V_0)$, such that
$-\Lambda_1-...-\Lambda_p$ is not an eigenvalue of
an element of ${\mathcal O}$ for any $0\le p\le N$.

Let $M_N$ be the set of collections $(a,b)$ of linear maps
$a_i: V_i\to V_{i+1}, b_i: V_{i+1}\to V_i$, $i=0,...,N$,
such that $b_ia_i-a_{i-1}b_{i-1}=\Lambda_i{\rm Id}_{V_i}$ for $i=1,...,N$,
and $b_0a_0\in {\mathcal O}$.

Let $x(a,b):=a_{N}b_{N}\in {\frak{gl}}(V_{N+1})$.
The group $G_N:=\prod_{i=0}^N GL(V_i)$ acts naturally on $M_N$
preserving the function $x(a,b)$.

\begin{lemma}\label{lina} (\cite{CB1})
(i) Let $(a,b)\in M_N$, and $x(a,b)=C$.
Then $C$ is conjugate to
$$
C_0\oplus \oplus_{i=1}^{N+1} (\Lambda_N+...+\Lambda_i){\rm Id}_{D_i-D_{i-1}}.
$$
where $C_0-(\Lambda_N+...+\Lambda_1)\in {\mathcal O}$
(here the subscripts $D_i-D_{i-1}$ denote matrix
sizes)\footnote{We agree that if $i=N+1$ then $\Lambda_N+...+\Lambda_i=0$}.

(ii) For any $C\in {\frak{gl}}(V_{N+1})$ as in (i),
there exists an element $(a,b)\in M_N$
such that $x(a,b)=C$. Moreover, any
two such elements are conjugate under $G_N$.
\end{lemma}

\begin{proof} The proof is by induction in $N$.
The base of induction ($N=0$) is easy.
Now assume that $N$ is arbitrary, and the
statement is known for $N-1$.

To prove (i), note that
by the induction assumption
$a_{N-1}b_{N-1}$ is conjugate to
$$
C_0'\oplus \oplus_{i=1}^{N} (\Lambda_{N-1}+...+\Lambda_i){\rm Id}_{D_i-D_{i-1}}.
$$
where $C_0'-(\Lambda_{N-1}+...+\Lambda_1)\in {\mathcal O}$.
But we have $a_{N-1}b_{N-1}=b_{N}a_{N}-\Lambda_N$,
So $b_Na_N$ is conjugate to
$$
C_0\oplus \oplus_{i=1}^{N} (\Lambda_{N}+...+\Lambda_i){\rm Id}_{D_i-D_{i-1}}.
$$
where $C_0-(\Lambda_{N}+...+\Lambda_1)\in {\mathcal O}$.
By our assumption, this implies that $b_Na_N$ is invertible, hence
$a_Nb_N$ is conjugate to the direct sum of $b_Na_N$ and the zero matrix
of size $D_{N+1}-D_N$, and (i) follows.

To prove (ii), pick $C'\in {\frak{gl}}(V_N)$ conjugate to
$$
C_0'\oplus \oplus_{i=1}^{N} (\Lambda_{N-1}+...+\Lambda_i){\rm Id}_{D_i-D_{i-1}}.
$$
where $C_0'-(\Lambda_{N-1}+...+\Lambda_1)\in {\mathcal O}$.
By the induction assumption, there exist unique up to conjugation
by $G_{N-1}$ operators $(a_0,...,a_{N-1},b_0,...,b_{N-1})\in M_{N-1}$
such that $a_{N-1}b_{N-1}=C'$. It remains to show that there exist
operators $a_N,b_N$ such that $a_Nb_N=C, b_Na_N=\Lambda_N+C'$,
and they are unique up to the action of the centralizer of $C'$ in $GL(V_N)$.
To do so, note that we must pick $a_N$ to be an isomorphism $V_N\to {\rm Im}C$,
which conjugates $C'+\Lambda_N$ to the restriction of $C$ to its image.
This can be done uniquely up to the action of the centralizer of $C'$, and
$b_N$ is uniquely determined by $a_N$. This implies the required statement.
\end{proof}

Now we prove the proposition. The existence of the map $\beta$ follows
from part (i) of Lemma \ref{lina}, by applying the Lemma
separately to each leg of $D$. Part (ii) of the lemma implies
that $\beta$ is bijective. But by Proposition \ref{smocl},
${\mathcal M}_{n,\gamma,\nu}$ is smooth, so $\beta$ is an
isomorphism, as desired.
\end{proof}

Now we will explicitly construct an isomorphism
$\Phi: {\mathcal R}_{n,\gamma,\nu}\to {\mathcal M}_{n,\gamma,\nu}$.
To do so, note that by Theorem \ref{mapro}, (iv),
every representation $V\in {\mathcal
R}_{n,\gamma,\nu}$ restricts to the regular representation
of ${\bold B}_{n,\ell}(\gamma_m,\nu)$ via the map $\eta_m$.
Thus the space $V':=V^{{\bold B}_{n-1,\ell}(\gamma_m,\nu)}$
has dimension $n\ell$. The operators
$Y_{n,1},...,Y_{n,m-1}$ commute with
${\bold B}_{n-1,\ell}(\gamma_m,\nu)$,
hence so does $Y_{n,m}-\nu(s_{n1}+s_{n2}+...+s_{n,n-1})$.
Thus, these elements define linear operators on $V'$.
Denote these operators by $x_1,...,x_m$, and set
$\Phi(V)=(x_1,...,x_m)$.

\begin{prop}\label{eig}
We have $(x_1,...,x_m)\in {\mathcal M}_{n,\gamma,\nu}$, so
$\Phi: {\mathcal R}_{n,\gamma,\nu}\to {\mathcal M}_{n,\gamma,\nu}$.
\end{prop}

\begin{proof} A simple deformation argument
from the case $\nu=0$ shows that the spectral decompositions of
$x_1,...,x_{m-1}$ are as required.
The fact that the spectral decomposition of $x_m$ is as required
follows from Lemma \ref{eigen}.
\end{proof}

\begin{prop}\label{comp}
$\Phi=\beta\circ \Phi_{EG}$.
\end{prop}

\begin{remark} Note that as a by-product
Proposition \ref{comp} gives another proof of
Proposition \ref{eig}.
\end{remark}

\begin{proof}
Let $W$ be an irreducible representation of
$\hh(\GG)$. In this case by Corollary \ref{cor20}
the corresponding representation $V$ of $B_n$
is $e_{i_0}^{\otimes n}W$. On the other hand,
the subspace $W^{S_{n-1}\ltimes \Gamma^{n-1}}$
considered in \cite{EG} equals\footnote{We note that
in \cite{EG} the subgroup $S_{n-1}\subset S_n$
is taken to be the stabilizer of $1$, while here it is taken to
be the stabilizer of $n$.} \linebreak
$e(n-1)(e_0^{\otimes (n-1)}\otimes 1)W$,
where $e(n-1)$ is the symmetrizer of $S_{n-1}$.

We need to construct an isomorphism
$\zeta: (e_{i_0})_nW^{S_{n-1}\ltimes \Gamma^{n-1}}\to V'$,
where $(e_{i_0})_n$ is the element $1^{\otimes (n-1)}\otimes e_{i_0}$.
This isomorphism is defined as follows.

Let $a_0,...,a_{\ell-1}$ be the opposites of the
edges of $Q$ belonging to the $m$-th leg ($a_{\ell-1}=h_m^*$).

Define $\zeta$ by the formula
$$
\zeta(w)=((a_{\ell-1}...a_1a_0)^{\otimes (n-1)}\otimes 1)w,
$$
where $(a_{\ell-1}...a_1a_0)^{\otimes (n-1)}\otimes 1)$ denotes an
element of $\mathsf{A}_n$. We note that
this element is well defined, because by Definition \ref{dea},
the elements $(b_1^*)_p$ ($b_1^*$ in the $p$-th tensor
component) and $(b_2^*)_q$ commute for any edges $b_1,b_2\in Q$
and $p\ne q$.

We claim that $\zeta(w)$ belongs to $V'$.
Indeed, $\zeta(w)$ is clearly invariant under $S_{n-1}$. Also,
for any $i<n$ we have $Y_{i,m}\zeta(w)=\gamma_{m\ell}\zeta(w)$.
To see this, recall that
$Y_{i,m}=(a_{\ell-1}a_{\ell-1}^*)_i+\xi_m+\mu_{i_0}/m$;
thus the statement follows from the relation (i) of Definition
\ref{dea} and formula (\ref{mugamma}).

It is easy to show that $\zeta$ is injective. Since $\zeta$ is a
morphism between spaces of the same dimension, it is an isomorphism.

It is now easy to check that if the spaces
$(e_{i_0})_nW^{S_{n-1}\ltimes \Gamma^{n-1}}$ and $V'$
are identified using $\zeta$, then the
quiver-theoretical data of \cite{EG}, Section 11 and the matrices
$x_1,...,x_m$ introduced above are related by the map $\beta$.
The proposition is proved.
\end{proof}

\begin{corollary}\label{is}
$\Phi$ is an isomorphism.
\end{corollary}

\begin{proof}
The corollary follows from Proposition \ref{comp} and
Proposition \ref{beta}.
\end{proof}

\section{Generalized double affine Hecke algebras}\label{s20}

\subsection{Generalized DAHA of rank 1}
First, let us recall the constructions in \cite{EOR,ER}.
Consider the group $G$ with generators $U_k$, $k=1, \ldots, m$,
and defining relations $$U_k^{d_k}=1,\quad k=1, \ldots, m,
\quad\textrm{and}\quad \prod_{k=1}^m U_k=1.$$
This group is a discrete group of motions
of the Euclidean plane, or Lobachevsky plane,
generated by rotations by the angles $2\pi/d_k$ around the
vertices of the $m$-gon with angles $\pi/d_k$, $k=1,...,m$.
Thus $G$ is a Euclidean crystallographic group
$\ZZ_\ell \ltimes \ZZ^2$ where $\ell=2,3,4,6$ for $D$ being
affine ($\widetilde D_4$, $\widetilde E_6$,
$\widetilde E_7$, $\widetilde E_8$ respectively),
and a hyperbolic motion group otherwise.

The generalized double affine Hecke algebra of rank 1 associated to
$D$ will be denoted by $H_1$. It is an algebra over
$\C[u^{\pm 1}]$, where $u$ denotes the
collection of variables
$$
u_{11}, \ldots, u_{1d_1},
\ldots, u_{m1}, \dots, u_{md_m}.
$$
The algebra $H_1$ is generated over $\C[u^{\pm 1}]$ by
elements $U_k$, $k=1,\ldots, m$, with defining relations
$$\prod_{j=1}^{d_k} (U_k-u_{kj})=0,\quad k=1,\ldots,m,\quad
\textrm{and}\quad \prod_{k=1}^m U_k=1.$$

The group algebra $\C[G]$ is isomorphic to the quotient of $H_1$ by the two-sided
ideal generated by $u_{kj}-e^{2\pi{\rm i}j/d_k}$, for all $k,j$.
In other words, $\C[G]$ is the specialization of $H_1$ at the
values $u_{kj}=e^{2\pi{\rm i}j/d_k}$.
Thus, $H_1$ is a deformation of $\C[G]$.

\subsection{Generalized DAHA of higher rank} We will now generalize the
definition of $H_1$ to the higher rank case.
Fix a positive integer $n>1$. Let $t$ be an additional invertible
variable.

\begin{definition}
The generalized double affine Hecke algebra $H_n$ of rank $n$ associated
to $D$ is the algebra generated over
$\C[u^{\pm 1}, t^{\pm 1}]$ by invertible elements
$$U_1, \ldots, U_m, T_1, \ldots, T_{n-1},$$ with defining relations
\begin{gather*}
(U_1\cdots U_m)(T_1T_2\cdots T_{n-2}T_{n-1}^2T_{n-2}\cdots
T_2T_1)=1,\\
T_iT_{i+1}T_i=T_{i+1}T_iT_{i+1},\quad i=1,\dots,n-2,\\
[T_i,T_j]=0,\quad |i-j|>1,\\
[U_j,T_i]=0,\quad i=2,\dots,n-1, \quad j=1,\dots,m,\\
[U_j,T_1U_jT_1]=0, \quad j=1,\dots,m,\\
[U_k,T_1^{-1}U_jT_1]=0, \quad 1\leq k<j\leq m,\\
\prod_{j=1}^{d_k}(U_k-u_{kj})=0,\quad k=1,\dots,m,\\
T_i-T_i^{-1}=t-t^{-1},\quad i=1,\dots,n-1.
\end{gather*}
\end{definition}

\begin{remark} The specialization of $H_n$ at the value
$t=1$ is the semidirect product  $\C[S_n]\ltimes H_1^{\ot n}$.
\end{remark}

\begin{remark} If we eliminate the last two groups of relations,
we get a presentation for the $n$-th braid group
$Br_{n,m}$ of $\C \Bbb P^1$ without $m$ points.
Thus the algebra $H_n$ can be viewed as a quotient
of the group algebra $\C[u^{\pm 1},t^{\pm 1}][Br_{n,m}]$ by the last two
groups of relations.
\end{remark}

\subsection{The case $D=\widetilde D_4$}
In the case when $D$ is of type $\widetilde D_4$,
the algebra $H_n$ is essentially the same as
the algebra $\mathcal{H}_n$ introduced by Sahi \cite[\S 3]{Sa},
which we now recall.

\begin{definition} $\mathcal{H}_n$ is the algebra
generated over $\C[t_0^{\pm 1},t_n^{\pm 1},u_0^{\pm 1},u_n^{\pm
1},t^{\pm 1},q^{\pm 1}]$
by elements $T_i^{\pm 1}$, $i=0,\dots,n$, and
elements $X_i^{\pm 1}$, $i=1,\dots,n$, subject to the relations
\begin{gather*}
T_0T_1T_0T_1=T_1T_0T_1T_0,\\
X_iX_j=X_jX_i,\ 1\le i<j\le n,\\
T_{n-1}T_nT_{n-1}T_n=T_nT_{n-1}T_nT_{n-1},\\
T_iT_{i+1}T_i=T_{i+1}T_iT_{i+1},\quad i=1,\dots,n-2,\\
[T_i,T_j]=0,\quad |i-j|>1,\\
T_i-T_i^{-1}=t-t^{-1},\quad i=1,\dots,n-1,\\
T_0-T_0^{-1}=t_0-t_0^{-1},\quad T_n-T_n^{-1}=t_n-t_n^{-1},\\
T_iX_j=X_jT_i,\mbox{ if }|i-j|>1,\mbox{ or }\mbox{ if } i=n
\mbox{ and } j=n-1,\\
T_iX_i=X_{i+1}T_i^{-1},\quad i=1,\dots,n-1,\\
T_n^{\vee}-(T_n^{\vee})^{-1}=u_n-u_n^{-1},\quad
\mbox{ where }T_n^{\vee}=X_n^{-1}T_n^{-1},\\
T_0^{\vee}-(T_0^{\vee})^{-1}=u_0-u_0^{-1},\quad
\mbox{ where }T_0^{\vee}=q^{-1}T_0^{-1}X_1.
\end{gather*}
\end{definition}

Let $H_n'$ be the specialization of $H_n$ defined by
$$
u_{11}=qt_0,\ u_{12}=-qt_0^{-1},\ u_{21}=u_0,\ u_{22}=-u_0^{-1},
$$
$$
u_{31}=u_n,\ u_{32}=-u_n^{-1},\ u_{41}=t_n,\ u_{42}=-t_n^{-1}.
$$
(This specialization is generic, in the sense that any set of
parameter values can be obtained from this one by the rescaling
transformations $u_{kj}\to u_{kj}w_k$.)

\begin{prop}\label{sahi}
There is an isomorphism $\phi: H_n'\to \mathcal{H}_n$,
given by
\begin{gather*}
\phi(U_1)=qT_0,\quad \phi(U_2)=T_0^{\vee},\quad \phi(U_3)=
ST_n^{\vee}S^{-1},\quad \phi(U_4)=ST_nS^{-1},\\
\phi(T_i)=T_i, \quad i=1, \ldots, n-1,
\end{gather*}
where $S=T_1T_2\dots T_{n-1}$.
\end{prop}

The proof of the proposition is by a direct computation.

\begin{remark}
The $T_i$'s from \cite{Sa} are denoted by $V_i$'s in \cite{St}.
Also, the following relations should be added in
\cite[Theorem 3.4]{St}:
\begin{gather*}
[V_0^{\vee},V_n]=[V_0,V_n^{\vee}]=0.
\end{gather*}
This is very minor misprint and it does not affect any other
results of the paper \cite{St} and the subsequent papers.
\end{remark}

\subsection{The flatness theorem}

\begin{conjecture}
The algebra $H_n$ is a free module over $\C[u^{\pm 1},t^{\pm 1}]$.
\end{conjecture}

In the case $n=1$, this conjecture is proved in \cite{EOR}.
Also, by Proposition \ref{sahi}
and the results of \cite{Sa}, the conjecture is true for
the affine diagram $\widetilde D_4$ for all $n$.

We can prove only a weaker version of this conjecture, which is
the following theorem. Let $\widehat H_n$ be the completion of
$H_n$ with respect to the ideal generated by $t-1$, and let
$\nu=-\frac{1}{\pi i}\log t$.

\begin{theorem}\label{def}
The algebra $\widehat H_n$ is a flat 1-parameter deformation
of the algebra $S_n\ltimes H_1^{\otimes n}$ with deformation
parameter $\nu$ (i.e., $\widehat H_n=S_n\ltimes H_1^{\otimes
n}[[\nu]]$ with deformed multiplication).
\end{theorem}

To prove this theorem, note first that by
a general deformation argument, it is sufficient to show
the following:

\begin{prop}\label{compl}
The completion $(\widehat H_n)_{u_0}$ of
$\widehat H_n$ at some point $u_0$ is a flat deformation
of the completion $(S_n\ltimes H_1^{\otimes n})_{u_0}$.
\end{prop}

We will give two proofs of this fact, using two different choices
of the point $u_0$.

{\bf First proof of Proposition \ref{compl}.} Let $Y$ be the
Euclidean or hyperbolic plane
which carries the action of the group $G$. Let $G_n=S_n\ltimes
G^n$. Then $G_n$ acts properly discontinuously on
$Y^n$, so following \cite{E}, we can define the Hecke algebra ${\mathcal
H}_\tau(Y^n,G_n)$ attached to the orbifold $Y^n/G_n$.
Using the explicit description of the braid group $Br_{n,m}$
of $Y^n/G_n$ given above, we find that the algebra
$(\widehat H_n)_{u_0}$ for $(u_0)_{kj}=e^{2\pi {\rm i}j/d_k}$
is a specialization of the algebra
${\mathcal H}_\tau(Y^n,G_n)$, $\tau=(\gamma,\nu)$. On the other hand
\footnote{Here is the place where
we use the fact that $D$ is not
of finite Dynkin type},
since $\pi_2(Y^n)=0$, one of the main results
of \cite{E} says that the algebra
$\mathcal H_\tau(Y^n,G_n)$ is flat over $\C[[\tau]]$.
This implies the theorem.

A second proof of Theorem \ref{def}, which does not use
the results of \cite{E}, will be given in the next section.
It is based on a variant of Knizhnik-Zamolodchikov equations.

\section{Knizhnik-Zamolodchikov equations}\label{s30}

\subsection{KZ equations}
Let $\alpha_1, \ldots, \alpha_m$ be distinct points in $\C$.
Consider the connection $\nabla$ on the trivial bundle
over $(\C\PP^1)^n$ with fiber $B_n$, defined by the
system of Knizhnik-Zamolodchikov (KZ) differential equations
\begin{equation}\label{kz}
\frac{\pa F}{\pa z_i}= A_i F, \quad i=1,\ldots, n,
\end{equation}
where
$$A_i := \sum_k \frac{Y_{i,k}}{z_i-\alpha_k}
        - \sum_{p\neq i}\frac{\nu s_{ip}}{z_i-z_p} \,.$$

\begin{lemma}
The connection $\nabla$ is flat.
\end{lemma}
\begin{proof}
To show that the curvature is zero, we have to
check that $$ \pa_i A_j-\pa_j A_i +[A_i,A_j]=0\,.$$
This follows from the following computations:
$$  \pa_j A_i = -\frac{\nu s_{ij}}{(z_i-z_j)^2} = \pa_i A_j \,,$$
\begin{align}
 [A_i,A_j]
=& \sum_{k,l}\frac{[Y_{i,k},Y_{j,l}]}{(z_i-\alpha_k)(z_j-\alpha_l)}
-\sum_{k, q\neq j}\frac{[Y_{i,k},\nu s_{jq}]}{(z_i-\alpha_k)(z_j-z_q)}
\nonumber\\
&\,-\sum_{k, p\neq i}\frac{[\nu s_{ip},Y_{j,k}]}{(z_i-z_p)(z_j-\alpha_k)}
+\sum_{p\neq i,q\neq j}\frac{[\nu s_{ip},\nu s_{jq}]}
{(z_i-z_p)(z_j-z_q)} \nonumber\\
=& \sum_k \frac{\nu(Y_{i,k}-Y_{j,k})s_{ij}}{(z_i-\alpha_k)(z_j-\alpha_k)}
-\sum_k \frac{\nu(Y_{i,k}-Y_{j,k})s_{ij}}{(z_i-\alpha_k)(z_j-z_i)}
\nonumber\\
&\,-\sum_k \frac{\nu(Y_{i,k}-Y_{j,k})s_{ij}}{(z_i-z_j)(z_j-\alpha_k)}
+\sum_{q\neq i,j}\frac{[\nu s_{ij},\nu s_{jq}]}{(z_i-z_j)(z_j-z_q)}
\nonumber\\
&\,+\sum_{p\neq i,j}\frac{[\nu s_{ip},\nu s_{ij}]}{(z_i-z_p)(z_j-z_i)}
+\sum_{p\neq i,j}\frac{[\nu s_{ip},\nu s_{jp}]}
{(z_i-z_p)(z_j-z_p)} \nonumber\\
=& 0 \,.\nonumber
\end{align}
\end{proof}

\subsection{The monodromy representation of the KZ equations}\label{monrep}

Taking quotient by the $S_n$-action, we get a flat connection,
which we will also denote by $\nabla$, on the configuration
space $\mathcal{C}onf_n(\C\PP^1\setminus
\lbrace{\alpha_1,...,\alpha_m\rbrace})$.
Note that we may replace
the trivial bundle over $(\C\PP^1)^n$ with fiber $B_n$ by the
trivial bundle whose fiber is a $B_n$-module $M$, and this also gives
a flat connection $\nabla_M$ on $\mathcal{C}onf_n$.

When the $B_n$-module
$M$ is finite dimensional, it acquires an action of the monodromy
operators. Namely, given a base point $\bold z_0\in
\mathcal{C}onf_n$, we can define the ${\rm End}M$-valued
solution $F_0$ of the KZ equations such that $F_0(\bold z_0)=1$.
Then, given $\sigma\in \pi_1(\mathcal{C}onf_n,\bold
z_0)$, we let $F_\sigma$ be the analytic continuation of $F_0$
along $\sigma$, and $L_\sigma\in {\rm End}M$ by
$F_\sigma=F_0L_\sigma$. Then the monodromy representation
$\rho: \pi_1(\mathcal{C}onf_n,\bold z_0)\to {\rm Aut}(M)$ is
defined\footnote{Our convention for the multiplication
of loops in $\pi_1$ is as follows:
to obtain $\sigma\sigma'$, first trace $\sigma'$,
then $\sigma$.} by the formula $\rho(\sigma)=L_\sigma$.

For convenience let us choose $\bold z_0=(z_{01},...,z_{0n})$ to be such that
$z_{0j},\alpha_p$ are real and $\alpha_1<...<\alpha_m<z_{01}<...<z_{0n}$.
In this case we can identify $\pi_1(\mathcal{C}onf_n,\bold z_0)$
with $Br_{n,m}$ as follows:
$T_i$ is the path in which the points $z_i,z_{i+1}$
move counterclockwise to exchange positions, and other points
don't move; $U_k$ is the path in which $z_1$ moves
counterclockwise around $\alpha_k$ (passing $\alpha_{k+1},...,\alpha_m$ from
below). Thus $\rho$ may be viewed as a representation of the
group $Br_{n,m}$ on $M$. 

Moreover, we claim that 
that the operators $T_i$ and
$U_k$ in this representation
satisfy the relations
$$
T_i-T_i^{-1}=t-t^{-1},
$$
and
$$
(U_k-u_{k1})...(U_k-u_{kd_k})=0,
$$
where
\begin{equation}\label{rattrig}
u_{kj}:=\exp(2\pi {\rm
i}\gamma_{kj}),\
t=e^{-\pi {\rm i}\nu}.
\end{equation}

Indeed, to prove the equation for $U_k$ it suffices to 
consider the KZ equation for the derivative with respect to 
$z_1$ with other variables fixed, and look at the 
eigenvalues of the residue of 
the connection at the point $z_1=\alpha_k$. On the other hand, to
prove the equation for $T_i$, transform the KZ equations
by the change of variables $z_{i,i+1}^+=\frac{z_i+z_{i+1}}{2}$, 
$z_{i,i+1}^-=\frac{z_{i+1}-z_i}{2}$, 
(leaving the variables other than $z_i,z_{i+1}$ unchanged); 
then the loop $T_i$ can be realized as a semicircle in which 
$z_{ii+1}^-$ goes from some (small) positive value $\zeta$ to
$-\zeta$ counterclockwise, and other variables (including
$z_{i,i+1}^+$) are unchanged. Looking at the eigenvalues of the residue of the
connection at $z_{i,i+1}^-=0$, we deduce the equation for $T_i$. 

Therefore, the monodromy representation of $Br_{n,m}$ on $M$ is in fact a
representation of the algebra $H_n$ with parameters $u,t$ as
above. Let us denote this representation of $H_n$ by $\F(M)$.
Thus we have obtained the following result.

\begin{prop} \label{prop30}
The monodromy of the KZ equations defines a functor
$\F:\Rep_f B_n \to \Rep_f H_n$ between the categories
of finite dimensional representations of $B_n$ and $H_n$, under
which the parameters $\gamma,\nu$ and $u,t$ are related as above.
\end{prop}

\begin{remark}
This functor, of course, depends on the choice of $\alpha_k$, but
only up to fractional-linear transformations.
\end{remark}

We note that this proposition allows us to construct a large supply of
finite dimensional representations of $H_n$ in the case when $D$
is affine. Indeed, a large supply of finite dimensional representations
for $\mathsf{A}_n$ and $\hh(\GG)$ (and hence for $B_n$) is
constructed in \cite{Ga,EM,M}, and we can apply the functor
$\F$ to these representations to obtain representations of
$H_n$.

{\bf Second proof of Proposition \ref{compl}.}

Let $\widetilde B_n$ be the formal completion of $B_n$ at the
point $\gamma=0,\nu=0$, and let $\widetilde H_n$
be the formal completion of $H_n$ at the ``unipotent point''
$t=1$, $u_{kj}=1$.

The monodromy of the KZ equation defines a morphism
$f: \widetilde H_n\to \widetilde B_n$, where parameters are
related as above (the Riemann-Hilbert homomorphism).
It is clear that this homomorphism is in fact an isomorphism
(since the relations of $B_n$ are infinitesimal versions of the
relations of $H_n$).
Thus, Proposition \ref{compl} follows from Corollary \ref{fre}.

\subsection{Cyclotomic Hecke algebras}

Let ${\bold H}_{n,\ell}$ be the Ariki-Koike cyclotomic Hecke
algebra (see e.g. \cite{Ma}). It is an algebra over $\C[v^{\pm 1},t^{\pm 1}]$
(where $v=(v_1,...,v_\ell)$) with generators $T_1,...,T_{n-1},U$
and defining relations
\begin{gather*}
T_iT_{i+1}T_i=T_{i+1}T_i T_{i+1}, \quad i=1,\dots, n-2,\\
[T_i, T_j]=0, \quad \mbox{ if }\quad |i-j|>1,\\
[U,T_j]=0,\quad j=2,\dots, n-1,\\
UT_1UT_1=T_1UT_1U,\\
\prod_{j=1}^{\ell}(U-v_j)=0,\\
 T_i-T_i^{-1}=t-t^{-1},\quad i=1,\dots,n-1.
\end{gather*}

Thus the algebra ${\bold H}_{n,\ell}$ is a quotient of the group
algebra $\C[v^{\pm 1},t^{\pm 1}][Br_{n,2}]$
(the fundamental group of the configuration space of $\Bbb C^*$)
by the polynomial relations for $U$ and $T_i$. It is known
(see e.g. \cite{Ma}) that ${\bold H}_{n,\ell}$ is a free module over
$\C[v^{\pm 1},t^{\pm 1}]$ of rank $n!\ell^n$, and thus
${\bold H}_{n,\ell}(v,t)$ is of dimension $n!\ell^n$
for all $v,t$ (i.e., it is a flat deformation of
the group algebra $\C[S_n\ltimes (\Bbb Z/\ell \Bbb Z)^n]$).

The algebra ${\bold H}_{n,\ell}(v,t)$ has a 1-dimensional
representation $\chi$ given by the formula
$\chi(T_i)=t$, $\chi(U)=v_\ell$. We call this representation the
trivial representation.
For generic parameters $v,t$,
the representation $\chi$ defines an
idempotent in ${\bold B}_{n,\ell}(\lambda,\nu)$.
We will denote this idempotent by ${\rm e}$.

For any representation
$V\in {\rm Rep} {\bold H}_{n,\ell}(v,t)$, we
denote by $V^{{\bold H}_{n,\ell}(v,t)}$ the space of
homomorphisms of representations $\chi\to V$.
Obviously, this space can be naturally regarded as a subspace of
$V$. In the generic case, this subspace is
equal to ${\rm e}V$.

Consider the subalgebra ${\bold H}_{n-1,\ell}$ of ${\bold H}_{n,\ell}$
generated by $U,T_1, T_2,...,T_{n-2}$.
For any representation $V$
of ${\bold H}_{n,\ell}(v,t)$, denote by $V'$ the space
$V^{{\bold H}_{n-1,\ell}(v,t)}$.

Let $X:=T_{n-1}...T_1UT_1...T_{n-1}\in {\bold
H}_{n,\ell}(v,t)$. In the braid group, the element $X$
corresponds to the point $z_n$ making a counterclockwise loop
around $0,z_1,...,z_{n-1}$; thus $X$ commutes with
${\bold H}_{n-1,\ell}(v,t)$. Therefore, $X$ acts
on the space $V'$ for any representation $V$
of ${\bold H}_{n,\ell}(v,t)$.

\begin{lemma}\label{eigen1} Let $V$ be the regular representation
of ${\bold H}_{n,\ell}(v,t)$. Then for generic parameters the
space $V'$ has dimension
$n\ell$, and the operator $X|_{V'}$ is conjugate to
$$
v_\ell t^{2T}
\oplus\mathop{{\rm diag}}(v_{1},\dots,v_{\ell-1})\otimes {\rm
Id}_n,
$$
where $T$ is as in Subsection \ref{dcha}.
\end{lemma}

\begin{proof} Let $V_0$ be the regular representation
of the degenerate cyclotomic Hecke algebra
${\bold B}_{n,\ell}(\lambda,\nu)$. Consider the following KZ differential
equations for a function $F(z_1,...,z_n)$ of complex variables
$z_1,...,z_n$ with values in $V_0$ (introduced by Cherednik \cite{Ch1}):
\begin{equation}\label{kz1}
\frac{\pa F}{\pa z_i}=\left(\frac{Y_i}{z_i}-
\sum_{p\ne i} \frac{\nu s_{ip}}{z_i-z_p}\right)F.
\end{equation}
(these are essentially equations (\ref{kz})
with $\alpha_m=0$, $\alpha_i=\infty$ for $i<m$).
Let $V$ be the monodromy representation
of this differential equation, with base point $\bold z_0=
(z_{01},...,z_{0n})$,
$0<z_{01}<z_{02}<...<z_{0n}$.
This is a representation of the braid group
$Br_{n,2}$ which obviously factors through the Hecke algebra
${\bold H}_{n,\ell}(v,t)$, where $t=e^{-\pi {\rm i}\nu}$
and $v_j=e^{2\pi{\rm i}\lambda_j}$. It is clear that generically
$V$ is the regular representation.

The space $V$ can be thought of as the space of local solutions
of the KZ equations around the base point.
The space $V^{{\bold H}_{n-1,\ell}(v,t)}$ can then be viewed as
the subspace of solutions $f$ of the form
$$
f=\prod_{j<i<n}(z_i-z_j)^{-\nu}\prod_{i<n}z_i^{\lambda_\ell}f_0,
$$
where $f_0$ analytically continues to a meromorphic function
in the region defined by the inequalities
$|z_n-z_{0n}|<\varepsilon$, $|z_i|<z_{0n}-\varepsilon$ for
$i<n$ (for some small $\varepsilon$). The operator
$X$ acts on the space of such solutions
by taking their monodromy around the loop $\sigma$
in which $z_n$ goes counterclockwise around $0,z_1,...,z_{n-1}$.
Tending $|z_i|$ to $0$, we find that the $n$-th KZ equation
tends to the equation $\frac{\pa F}{\pa z_n}=
\frac{Y_n-\nu(s_{n1}+...+s_{n,n-1})}{z_n}F$. Therefore, the
monodromy around $\sigma$ on $V'$ is conjugate
to $e^{2\pi {\rm i}(Y_n-\nu(s_{n1}+...+s_{n,n-1}))}|_{V_0'}$.
Thus the required statement follows from
Lemma \ref{eigen}.
\end{proof}

Here is another, purely algebraic proof of Lemma \ref{eigen1}.

\begin{proof}
Let $\pi\in {\bold
H}_{n-1,\ell}(v,t)$ be the
idempotent of the trivial representation (it exists since
generically the algebra ${\bold H}_{n-1,\ell}(v,t)$ is
semisimple); we have $V'=\pi V$.

Let $P_k(x)=\prod_{j=1,j\ne
k}^{\ell}\frac{(x-v_{j})}{(v_{k}-v_{j})}$.
Let $U_n=T_{n-1}^{-1}...T_1^{-1}UT_1...T_{n-1}$.
We have $[\pi,U_n]=0$, and hence
$V'=\oplus_{k=1}^\ell V_k'$, where $V_k':=P_k(U_n)V'$.
Note that $V_k'=y_kV'$, where $y_k=\pi
T_{n-1}^{-1}...T_1^{-1}P_k(U)$, and that $\dim V_k'=n$ for all
$k$.

\begin{lemma}\label{ei} $X|_{V_k'}=v_k{\rm Id}$ for any $k\ne \ell$.
\end{lemma}

\begin{proof}
It suffices to show that
$$Xy_{k}=v_{k}y_{k},
\quad k=1,\dots,\ell-1.$$
Since $X\pi=\pi X$, we find
$$
Xy_{k}=v_k\pi T_{n-1}...T_2T_1P_k(U)=
v_k\pi T_{n-1}...T_2T_1^{-1}P_k(U)+v_k(t-t^{-1})\pi T_{n-1}...T_2P_k(U).
$$
Since $P_k(U)$ commutes with $T_{n-1},...,T_2$ and $\pi
P_k(U)=0$, the last summand is zero. Thus we have
$$
Xy_{k}=v_k\pi T_{n-1}...T_2T_1^{-1}P_k(U)=
$$
$$
v_k\pi T_{n-1}...T_2^{-1}T_1^{-1}P_k(U)+
v_k(t-t^{-1})\pi T_{n-1}...T_3T_1^{-1}P_k(U),
$$
and again the last term is zero.
Continuing in this way, we will find that
$$
Xy_{k}=v_k\pi T_{n-1}^{-1}...T_2^{-1}T_1^{-1}P_k(U)=v_ky_{k}.
$$
\end{proof}

Now let ${\bold H}_n(t)$ be the usual Hecke algebra of type
$A_{n-1}$, and consider the homomorphism
$\theta: {\bold H}_{n,\ell}(v,t)\to {\bold H}_n(t)$ given by
$T_i\to T_i$, $U\to v_\ell$. We have $\theta(X)=v_\ell
T_{n-1}...T_2T_1^2T_2...T_{n-1}$. Also, $\theta(V_k')=0$ for $k\ne
\ell$, while $\theta|_{V_\ell'}$ is injective, and its image $J$ is
the space of all elements $y$ in ${\bold H}_n(t)$ such that
$T_iy=ty$ for $i\ne n-1$. By the result of Subsection 4.3 in
\cite{Ob1}, the operator $T_{n-1}...T_2T_1^2T_2...T_{n-1}|_J$
is conjugate to $t^{2T}$. This statement, together with Lemma
\ref{ei} implies Lemma \ref{eigen1}.
\end{proof}

\section{The affine case}

\subsection{The algebra $H_n$ for affine $D$}

 From now on let us consider the affine case, i.e.
$D=\widetilde D_4,\widetilde E_6,\widetilde E_7,\widetilde E_8$.
We keep the notation of subsection \ref{affi}.

It is natural to expect that in this case the algebra $H_n(u,t)$
has properties similar to those of the algebra $B_n(\gamma,\nu)$
stated in Theorem \ref{mapro}. Unfortunately, we are unable to
establish any of these properties, and
we are going to state them as conjectures.

Set
$$
q=q(u):=\prod_{k,j}u_{kj}^{-\ell/d_k}.
$$

We have a homomorphism
$\eta_m: {\bold H}_{n,\ell}(u_m,t)\to H_n(u,t)$, where
$u_m:=(u_{mj})$, given by
the formulas $\eta_m(T_i)=T_i$, $\eta_m(U)=U_m$.

\begin{conjecture} (i) The Gelfand-Kirillov dimension
of $H_n(u,t)$ is $2n$.

(ii) The algebra $H_n(u,t)$ is PI if and only
if $q$ is a root of unity. More precisely, it is PI of degree
$n!(N\ell)^n$ if $q$ is a root of unity of order $N$.

(iii) If $q$ is a root of unity, then $H_n(u,t)$ is finitely
generated over its center $Z(H_n(u,t))$.

(iv) If $q$ is a root of unity and otherwise $(u,t)$ are generic
then $H_n(u,t)$ is an Azumaya algebra, and $S(u,t):={\rm Spec}(Z(H_n(u,t)))$
is a smooth affine algebraic variety of dimension $2n$.

(v) If $q=1$ then the map $Z(H_n(u,t))\to {\rm e}H_n(u,t){\rm e}$
given by $z\mapsto z{\rm e}$ is an isomorphism. In particular,
${\rm e}H_n(u,t){\rm e}$ is a commutative algebra.

(vi) If $q=1$ and otherwise $(u,t)$ are generic then
every irreducible representation of $H_n(u,t)$
restricts (via the map $\eta_m$) to the regular representation of
$\bold H_{n,\ell}(u_m,t)$.
\end{conjecture}

\begin{remark}
For $n=1$, this conjecture follows from
the paper \cite{EOR}. Also, for $D=\widetilde D_4$,
because of Proposition \ref{sahi}
this conjecture can be attacked using the methods of \cite{Sa} and \cite{Ob1}
(for example, parts (i)-(iii) and the second statement of (v)
follow rather easily from \cite{Sa}); this will be done in a
subsequent paper. Finally, using the Riemann-Hilbert
homomorphism discussed in Section 4, it can be shown that
parts (i)-(iii) and (v) of the conjecture hold for the completed
algebra $\widehat H_n$. On the other hand, for the algebra $H_n$
of types $\widetilde E_l$, $l=6,7,8$,
it is unclear to us how to attack any of the above
questions (basically, because we don't know how to construct a
basis or at least a well behaved filtration of $H_n$, similar to
those used in \cite{EOR} for $H_1$).
\end{remark}

\subsection{Representations of $H_n(u,t)$ for $q=1$}

Assume that $q=1$ and otherwise $(u,t)$ are generic.
In this case, the algebra $H_n(u,t)$ has a $2n$-parameter family
of representations of dimension $n!\ell^n$, which are constructed
as follows. Let $\gamma,\nu$ satisfy equations (\ref{rattrig}), and
$\hbar=0$ (this is possible since $q=e^{-\hbar}$). In this case, by
Theorem \ref{mapro}, $B_n(\gamma,\nu)$ is an Azumaya algebra, so
all irreducible representations of
$B_n(\gamma,\nu)$ have dimension $n!\ell^n$, and
are parametrized by a smooth connected $2n$-dimensional algebraic variety
${\mathcal R}_{n,\gamma,\nu}$. Thus for any $M\in {\mathcal
R}_{n,\gamma,\nu}$, we can define a representation $\F(M)$ of
$H_n(u,t)$, of dimension $n!\ell^n$ (see Subsection \ref{monrep}).

\begin{prop}\label{irrres}
For generic $M\in {\mathcal R}_{n,\gamma,\nu}$,
$\F(M)$ is irreducible, and $\eta_m^*\F(M)$ is the regular
representation of ${\bold H}_{n,\ell}(u_m,t)$.
\end{prop}

\begin{proof} In the case $\nu=0$ the statement reduces to the
rank 1 case and hence follows from the results of \cite{EOR}.
Therefore, the statement holds for generic parameters and generic
$M$.
\end{proof}

Let ${\Bbb  R}_{n,u,t}$ be the set of equivalence classes of
irreducible representations of $H_n(u,t)$ which restrict
(via the map $\eta_m$) to the regular representation of
${\bold H}_n(u_m,t)$. This is an affine algebraic variety.
By Proposition \ref{irrres}, for generic $M$ as above,
$\F(M)\in {\Bbb  R}_{n,u,t}$.

\begin{remark}
As we mentioned in the previous subsection,
we conjecture that all irreducible representations
of $H_n(u,t)$ (for $q=1$ and otherwise generic $u,t$)
restrict to the regular representation of
${\bold H}_n(u_m,t)$ and thus belong to ${\Bbb  R}_{n,u,t}$.
\end{remark}

We now want to parametrize irreducible representations
of $H_n(u,t)$, by constructing a map $\Phi: {\Bbb  R}_{n,u,t}\to
{\Bbb  M}_{n,u,t}$ of ${\Bbb  R}_{n,u,t}$ into some explicitly
described algebraic variety ${\Bbb  M}_{n,u,t}$, similarly
to the map $\Phi$ for $B_n(\gamma,\nu)$ discussed in Subsection \ref{rep0}.

The variety ${\Bbb  M}_{n,u,t}$ is defined as follows.

\begin{definition} ${\Bbb  M}_{n,u,t}$
is the variety of conjugacy classes of $m$-tuples \linebreak
$(X_1,\dots,X_m)\in
GL_{n\ell}(\C)^{m}$ satisfying the following equations:
\begin{gather}\label{produ}
X_1X_2\dots X_m=1,\\
\label{U_12..}
X_k\sim \mathop{{\rm diag}}(u_{k1},\dots,u_{k,d_k})\otimes
{\rm Id}_{n\ell/d_k},\quad k=1,\dots,m-1,\\
\label{X_m}
X_m\sim u_{m\ell}t^{2T}
\oplus\mathop{{\rm diag}}(u_{m1},\dots,u_{m,\ell-1})\otimes {\rm Id}_n.
\end{gather}
That is, ${\Bbb  M}_{n,u,t}$ is the categorical quotient
of the variety $\widetilde {\Bbb  M}_{n,u,t}$
of $m$-tuples as above by the action of the group
$PGL_{n\ell}(\Bbb C)$.
\end{definition}

\begin{remark}
Thus, ${\Bbb  M}_{n,u,t}$ is defined as the variety
of solutions of an appropriate multiplicative Deligne-Simpson problem.
\end{remark}

\begin{prop}\label{smo} For generic parameters ${\Bbb  M}_{n,u,t}$ is a smooth
variety, of dimension $2n$.
\end{prop}

\begin{proof} The proof is standard (see also \cite{CBS}).
First of all, for generic parameters, any matrices
$X_1,...,X_m$ satisfying equations (\ref{produ},\ref{U_12..},\ref{X_m})
form an irreducible family. Indeed it is easy to see by
computing determinants of both sides of (\ref{produ}) using
equations (\ref{U_12..},\ref{X_m}) that
the only nonzero invariant subspace for $X_1,...,X_m$ is the
whole space. This implies that
the group $PGL_{n\ell}(\Bbb C)$ acts freely on
the variety $\widetilde {\Bbb  M}_{n,u,t}$.

It remains to show that the variety $\widetilde {\Bbb  M}_{n,u,t}$
is smooth, of dimension $2n+n^2\ell^2-1$. To do so,
let $C_1,...,C_m$ denote the conjugacy classes
of $X_1,...,X_m$. We have a map $\mu: C_1\times...\times C_m\to
SL_{n\ell}(\Bbb C)$ given by $(X_1,...,X_m)\to X_1...X_m$,
and $\widetilde {\Bbb  M}_{n,u,t}=\mu^{-1}(1)$. We have
$$
\dim C_k=n^2\ell^2(1-1/d_k),\ k=1,...,m-1;\quad
\dim C_m=2n-2+n^2\ell^2(1-1/d_m).
$$
Since for affine $D$, $\sum_k (1-1/d_k)=2$, we have
$$
\dim(C_1\times...\times C_m)=2n-2+2n^2\ell^2.
$$
Thus, to prove the proposition, it suffices to
show that $1$ is a regular value for the map $\mu$, i.e. that
for any $\bold X=(X_1,...,X_m)\in \widetilde {\Bbb  M}_{n,u,t}$, the
differential $d\mu_{\bold X}$ is surjective.

A tangent vector to $\bold X$ in $C_1\times...\times C_m$
is of the form $([P_1,X_1],...,[P_m,X_m])$, where $P_k$ are some
matrices. We have
$$
d\mu_{\bold X}([P_1,X_1],...,[P_m,X_m])=
\sum_{k=1}^mX_1...X_{k-1}[P_k,X_k]X_{k+1}...X_m=
$$
$$
\sum_{k=0}^{m}X_1...X_{k}(P_{k+1}-P_{k})X_{k+1}...X_m
$$
(where we agree that $P_0=P_{m+1}=0$).
Let $Q_k:=P_{k+1}-P_{k}$, $k=1,...,m$
(they can be arbitrary matrices). Then we get
$$
d\mu_{\bold X}([P_1,X_1],...,[P_m,X_m])=
\sum_{k=1}^m [X_1...X_{k},Q_k]X_{k+1}...X_m=
$$
$$
\sum_{k-=1}^m {\rm Ad}(X_1...X_k)(Q_k)-Q_k
$$
(as $X_1...X_m=1$). Now, since $X_1,...,X_m$ is an irreducible family,
we have $\cap_{k=1}^m {\rm Ker}({\rm Ad}(X_1...X_k)-1)=\Bbb C$, and
hence dually $\sum_{k=1}^m {\rm Im}({\rm Ad}(X_1...X_k)-1)={\frak
{sl}}_{n\ell}(\Bbb C)$. Thus, $d\mu_{\bold X}$ is surjective and we
are done.
\end{proof}

\begin{remark}
We have not shown that the variety
${\Bbb  M}_{n,u,t}$ is nonempty. This will follow from the
existence of the map
$\Phi$ defined below, and also follows from
the results of \cite{CBS}.
\end{remark}

Let us explain the construction of the map $\Phi$: ${\Bbb
R}_{n,u,t}\to {\Bbb M}_{n,u,t}$. Let $V\in {\Bbb
R}_{n,u,t}$. Using the map $\eta_m$ we can regard $V$ as a
representation of ${\bold H}_{n,\ell}(u_m,t)$, which is
isomorphic to the regular representation. It
is easy to see that the elements $\widetilde U_i:=
T_{n-1}...T_1U_iT_1^{-1}...T_{n-1}^{-1}$, $i=1,\dots,m-1$, commute
with ${\bold H}_{n-1,\ell}(u_m,t)$. The same is true about
$\widetilde{U}_m:=T_{n-1}\dots T_1U_mT_1\dots T_{n-1}$.
Thus the operators $\widetilde {U}_1,\dots,\widetilde {U}_{m-1},
\widetilde{U}_m$ preserve the space of
$V':=V^{{\bold H}_{n-1,\ell}(u_m,t)}$. We have
$\dim V'=n\ell$. We define the map $\Phi$ by the formula:
\begin{equation*}
\Phi(V)=(\widetilde{U}_1|_{V'},\dots,
\widetilde{U}_{m-1}|_{V'},\widetilde{U}_m|_{V'}).
\end{equation*}

\begin{prop} We have $\Phi: {\Bbb R}_{n,u,t}\to {\Bbb M}_{n,u,t}$.
\end{prop}

\begin{proof} An easy deformation argument from
the group algebra case shows that equation (\ref{U_12..}) is satisfied
for $X_k=\widetilde{U}_k|_{V'}$, $k<m$. Equation (\ref{produ}) is also clearly
satisfied. Finally, the fact that
equation (\ref{X_m}) holds follows from Lemma \ref{eigen1}.
\end{proof}

\begin{conjecture}
$\Phi$ is an isomorphism of algebraic varieties.
\end{conjecture}

Recall now that there is a Riemann-Hilbert map
between the spaces of solutions of the additive and the
multiplicative Deligne-Simpson problems: $RH: {\mathcal
M}_{n,\gamma,\nu}\to {\Bbb M}_{n,u,t}$, defined as follows.

Given $\bold x:=(x_1,...,x_m)\in {\mathcal
M}_{n,\gamma,\nu}$, consider the Fuchsian differential equation
\begin{equation}\label{rh}
\frac{dF}{dz}=\sum_{k=1}^m\frac{x_kF}{z-\alpha_k}.
\end{equation}

Assume that $z_0\in \Bbb R$ is a base point, and
$\alpha_1<...<\alpha_m<z_0$.
Then $RH(\bold x)=\bold X:=(X_1,...,X_m)$, where $X_k$
is the monodromy matrix of this differential equation
around the loop, in which
$z$ goes counterclockwise around $\alpha_k$ passing
$\alpha_{k+1},..,\alpha_m$ from below.
This map, of course, depends on the choice of $\alpha_k$, but
only up to fractional-linear transformations.

\begin{prop} One has $\Phi\circ {\mathcal F}=RH\circ \Phi$.
\end{prop}

\begin{proof} The proof is similar to the
first proof of Lemma \ref{eigen1}.
If $V\in {\mathcal R}_{n,\gamma,\nu}$, then
$\Phi\circ \F(V)$ is the collection of operators $\widetilde U_k$
on the subspace $\F(V)'$ of $\F(V)$. This subspace
can be viewed as the space of solutions of the KZ equations
which become single-valued near $z_i=\alpha_m$ and $z_i=z_j$ ($i,j<n$)
after division by
$\prod_{j<i<n}(z_i-z_j)^{-\nu}\prod_{i<n}(z_i-\alpha_m)^{\gamma_{m\ell}}$,
and $\widetilde{U}_k$ are the monodromy operators for
such solutions around the loops $\sigma_k$, in which
$z_n$ goes counterclockwise around $\alpha_k$ passing
$\alpha_{k+1},..,\alpha_m,z_1,...,z_{n-1}$ from below for $k<m$,
and goes around $\alpha_m,z_1,...,z_{m-1}$ for $k=m$.
To compute the spectral type of $\widetilde U_k$, we may send
$z_i$ with $i<n$ to zero. In this case, the $n$-th KZ equation
tends to equation (\ref{rh}). This implies the required
statement.
\end{proof}

In conclusion we would like to discuss the dependence of the map
$RH$ on the parameters $\alpha_k$. In the $\widetilde E_l$-cases,
there is no such dependence, since there are only three
parameters $\alpha_k$, and all collections of them are
projectively equivalent. On the other hand, in the $\widetilde
D_4$ case, we have an essential parameter, which is the
cross-ratio $\kappa$ of $\alpha_1,\alpha_2,\alpha_3,\alpha_4$.
Thus we have a 1-parameter family of holomorphic maps
$RH_\kappa: {\mathcal M}_{n,\gamma,\nu}\to {\Bbb M}_{n,u,t}$.
If we fix $\bold X\in {\Bbb M}_{n,u,t}$, we can (locally) implicitly solve
for $\bold x\in {\mathcal M}_{n,\gamma,\nu}$ such that
$RH_\kappa(\bold x)=\bold X$. This gives a function
$\bold x=\bold x(\kappa,\bold X)$, which defines a flow
on the 2n-dimensional complex manifold
${\mathcal M}_{n,\gamma,\nu}$. In the case $n=1$, this is the
Painlev\'e VI flow; so in general this flow should be regarded as
a higher rank version of Painlev\'e VI. Note that
the higher rank Painlev\'e VI flow has an additional parameter
$\nu$, so it has 5 parameters, rather than 4 for the usual
Painlev\'e VI; if $\nu=0$, the higher rank Painlev\'e
VI flow decouples into a (symmetric) product of $n$ copies of the usual
Painlev\'e VI flows. It would be interesting to write this differential
equation explicitly using an appropriate coordinate system
on ${\mathcal M}_{n,\gamma,\nu}$.
\setcounter{equation}{0}

\end{document}